\documentclass[12pt]{amsart}
\usepackage{pb-diagram,lamsarrow,pb-lams,amssymb}
\title[Loop space Bockstein spectral sequence]
  {Algebraic structure of the loop space \\
   Bockstein spectral sequence}
\author{Jonathan A. Scott}
\address{Department of Mathematics, University of Toronto,
  M5S 3G3 Canada}
\date{\today}
\email{scott@math.toronto.edu}
\newtheorem{thm}{Theorem}
\newtheorem{prop}[thm]{Proposition}
\newtheorem{lem}[thm]{Lemma}

\theoremstyle{definition}
\newtheorem{defn}{Definition}

\theoremstyle{remark}

\newtheorem{ex}{Example}

\newcommand{\tensor}{\otimes}

\newcommand{\dirsum}{\oplus}
\newcommand{\iso}{\cong}
\newcommand{\quism}{\stackrel{\simeq}{\rightarrow}}
\newcommand{\dual}{\sharp}
\newcommand{\zlocp}{\mathbf{Z}_{(p)}}
\newcommand{\fp}{\mathbf{F}_{\!p}}
\newcommand{\Q}{\mathbf{Q}}
\newcommand{\Z}{\mathbf{Z}}

\newcommand{\del}{\partial}
\newcommand{\ab}{L_{\mbox{\scriptsize{ab}}}}

\newcommand{\Hom}{\operatorname{Hom}}

\newcommand{\tor}{\operatorname{Tor}}

\newcommand{\image}{\operatorname{im}}

\newcommand{\dgh}{\textsc{dgh}}
\newcommand{\dgm}{\textsc{dgm}}
\newcommand{\dgl}{\textsc{dgl}}
\newcommand{\dga}{\textsc{dga}}
\newcommand{\dgg}{$\Gamma$-\textsc{dga}}
\newcommand{\dgc}{\textsc{dgc}}
\newcommand{\bss}{\textsc{bss}}

\begin{document}
\begin{abstract}
Let $X$ be a finite, $n$-dimensional, $r$-connected CW complex.
We prove the following theorem:

If $p \geq n/r$ is an odd prime, 
then the loop space homology
Bockstein spectral sequence modulo $p$ is a spectral sequence of
universal enveloping algebras over differential graded Lie algebras.
\end{abstract}

\maketitle
\section*{Introduction}

Let $\Omega X$ be the Moore loop space on a pointed topological
space $X$.
If $R \subseteq \Q$ is a principal ideal domain, 
then $H_{*}(\Omega X;R)$ has a natural Hopf algebra structure
via composition of loops, as long as there is no torsion.
The submodule $P \subset H_{*}(\Omega X; R)$ of primitive elements 
is a graded Lie subalgebra;
in~\cite{milnor-moore:65}, Milnor and Moore showed that if 
$R = \Q$ and $X$ is simply connected then
$H_{*}(\Omega X;\Q)$ is the universal enveloping algebra of $P$.
In~\cite{halperin:92}, Halperin established the same conclusion for
$R \subset \Q$ when $X$ is a finite, 
simply-connected CW complex, provided that 
$H_{*}(\Omega X;R)$ is torsion-free 
and the least non-invertible prime in $R$ is
sufficiently large.

In the presence of torsion, the loop space
homology algebra does not have a natural Hopf algebra
structure.
However, in~\cite{browder:61} Browder showed that the 
Bockstein spectral sequence
\[
  H_{*}(\Omega X;\fp) \Rightarrow
    \left(H_{*}(\Omega X;\Z)/\mbox{torsion}\right)\tensor\fp
\]
is a spectral sequence of Hopf algebras.
Halperin also proved in~\cite{halperin:92} that for large enough primes,
$H_{*}(\Omega X;\fp)$ is the universal enveloping algebra
of a graded Lie algebra.
The present article establishes this for every term in
the Bockstein spectral sequence.
\begin{thm}\label{thm:loopbss}
Let $X$ be a finite, $n$-dimensional, $q$-connected 
CW complex ($q \geq 1$).
If $p$ is an odd prime and $p \geq n/q$, then each term in the mod $p$
homology Bockstein spectral sequence for $\Omega X$ is
the universal enveloping algebra of a differential graded Lie
algebra $(L^{r},\beta^{r})$.
\end{thm}

In~\cite{anick:89}, under the hypotheses of 
Theorem~\ref{thm:loopbss}, Anick associates to $X$ 
a differential graded Lie algebra $L_{X}$ over $\zlocp$ and a natural
quasi-isomorphism $UL_{X} \quism C_{*}(\Omega X; \zlocp)$
of graded Hopf algebras.
The inclusion $\iota_{X}:L_{X} \hookrightarrow UL_{X}$
then induces a transformation of Bockstein spectral sequences
\(
  E^{r}(\iota_{X}):E^{r}(L_{X}) \rightarrow E^{r}(\Omega X).
\)
\begin{thm}\label{thm:transf}
The image of each $E^{r}(\iota_{X})$ is contained in $L^{r}$.
\end{thm}

Theorems~\ref{thm:loopbss} and~\ref{thm:transf} follow from
the work of Anick in~\cite{anick:89} and the following:
\begin{thm}\label{thm:main}
Let $(L,\del)$ be a differential graded Lie algebra
over $\zlocp$ which is connected,
free as a graded module, and of finite type.
The mod $p$ homology Bockstein spectral sequence of $U(L,\del)$
is a sequence of universal enveloping algebras,
\(
  E^{r}(UL) = U(L^{r},\beta^{r}).
\)
Furthermore, if $\iota:L\hookrightarrow UL$ is the inclusion,
then the image of $E^{r}(\iota)$ is contained in $L^{r}$.
\end{thm}

The proof of Theorem~\ref{thm:main} utilizes in a fundamental
way the divided powers structure of the dual of a universal
enveloping algebra.

The structure of the article is as follows.

\emph{Section~\ref{s:prelims}.} 
Notation and review of graded Lie algebras,
divided powers algebras, Bockstein spectral sequences,
acyclic closures and minimal models.

\emph{Section~\ref{s:lie-morph}.}
In~\cite{halperin:92}, Halperin showed that for a differential
graded Lie algebra $(L,\del)$ over $\fp$, 
$H(UL)=UE$ for a graded Lie algebra $E$.
We show that the inclusion 
$\iota:(L,\del)\hookrightarrow U(L,\del)$
satisfies $\image H(\iota) \subset E$.

\emph{Section~\ref{s:bss}.} 
Proof of Theorem~\ref{thm:main}.

\emph{Section~\ref{s:morphisms}.}
We show that a Hopf algebra morphism 
$UL_{1}\rightarrow UL_{2}$ is of the form $U(\varphi)$
if and only if its dual respects divided powers.
 
\emph{Section~\ref{s:examples}.}
Two examples.
The first gives a differential graded Lie algebra
whose 
Bockstein spectral sequence collapses after the first term,
while the spectral sequence of its universal enveloping algebra never
does.
The second shows that the sequence of Lie algebras 
given by Theorem~\ref{thm:main} is not natural.

\subsection*{Acknowledgments}

This paper is the result of work for my Ph.D. thesis.
I would like to take this opportunity to thank my Ph.D.
advisor, Steve Halperin, who suggested the problem of
studying torsion in loop space homology, and who patiently
and consistently provided encouragement, deep insight, and support.


\section{Preliminaries}\label{s:prelims}

Let $R$ be a commutative ring in which $2$ is invertible.
All objects are graded by the integers unless
otherwise stated.
Fix an odd prime $p$.
The ring of integers localized at $p$ is denoted
$\zlocp$ while the prime field is denoted $\fp$.
Differential graded modules, algebras, coalgebras, and Hopf algebras
are shortened to \dgm, \dga, \dgc, and \dgh, respectively;
a comprehensive treatment of these objects is given in
~\cite{f-h-t:95}.

\subsection{Graded modules}
Let $M$ be a graded module over $R$.
If $x \in M_{k}$ then we say that $x$ has degree $k$, and write
$|x| = k$.
A free graded module $M$ is of \emph{finite type} if each
$M_{k}$ is of finite rank.
We raise and lower degrees by the convention $M^{k}=M_{-k}$.
We denote by $sM$ the suspension of $M$: $(sM)_{i} = M_{i-1}$.
The dual of $M$ is the graded module 
$M^{\dual} = \Hom(M,R)$.
If $M$ is finite type and $N = (sM)^{\dual}$, then $M = (sN)^{\dual}$ via
$x(sf) = - f(sx)$, for $x \in M$, $f \in N$.

If $V$ is a graded module over $R$, then we denote by $TV$ and
$\Lambda V$ the tensor algebra and free commutative algebra
on $V$, respectively.
The tensor coalgebra on $V$ is denoted by $T_{C}V$.
The shuffle product (\cite{halperin:92}, Appendix) makes
$T_{C}V$ into a graded commutative (not cocommutative)
Hopf algebra.
Note that $TV = \dirsum_{k\geq 0}T^{k}V$,
$\Lambda V = \dirsum_{k\geq 0}\Lambda^{k}V$
and $T_{C}V = \dirsum_{k\geq 0}T_{C}^{k}V$, with
$T^{k}V$, $\Lambda^{k}V$, and $T_{C}^{k}V$ consisting of words in 
$V$ of length $k$.
Elements of $T_{C}^{k}V$ are denoted
$[v_{1}|\cdots|v_{k}]$.

The symmetric group $S_{k}$ acts on $T^{k}V$ via
$\sigma\cdot(x_{1} \tensor \cdots \tensor x_{k})
  = \pm\, x_{\sigma(1)} \tensor \cdots \tensor x_{\sigma(k)}$,
where the sign is determined by the rule
$x \tensor y \mapsto (-1)^{|x||y|}y \tensor x$.
\subsection{Graded Lie algebras}

A \emph{graded Lie algebra} is a graded $R$-module 
$L = \dirsum_{k\geq 0}L_{k}$ along with a degree-zero linear map
$[\,,]:L \tensor L \rightarrow L$, called the \emph{Lie bracket},
satisfying graded anti-commutativity, the graded Jacobi identity,
and the further condition $[x,[x,x]]=0$ if 
$x \in L_{\mbox{\scriptsize{odd}}}$.

For example, any non-negatively graded 
associative algebra $A$ is a graded Lie
algebra via the graded commutator bracket 
$[a,b] = ab - (-1)^{|a||b|}ba$, for $a,b \in A$.

A graded Lie algebra is \emph{connected}
if it is concentrated in strictly positive degrees.

The \emph{graded abelian Lie algebra} on
$\{x_{j}\}$, denoted
$\ab(x_{j})$,
is the free graded module on
the basis $\{x_{j}\}$, with the trivial Lie bracket.

Let $L$ be a graded Lie algebra, and denote by $L^{\flat}$ the 
underlying graded module.
The \emph{universal enveloping algebra} of $L$ is the
associative algebra $UL = (TL^{\flat})/I$, where
$I$ is the ideal generated by elements of the form
$x \tensor y - (-1)^{|x||y|}y \tensor x - [x,y]$,
for $x,y\in L$.
$UL$ has the natural structure of a graded Hopf algebra;
the comultiplication is defined by
declaring the elements of $L$ to be primitive and then using
the universal property.

A \emph{Lie derivation} on a graded Lie algebra $L$ is a linear operator 
$\theta$ on $L$ of degree $k$ such that for $x,y \in L$,
\(
  \theta([x,y]) = [\theta(x),y] + (-1)^{k|x|}[x,\theta(y)].
\)
A \emph{differential graded Lie algebra} (\dgl\ for short)
is a pair $(L,\del)$, where $L$ is a graded Lie algebra,
and $\del$ is a Lie derivation on $L$ of degree $-1$ satisfying
$\del\del = 0$.
If $(L,\del)$ is a \dgl, then $\del$ extends to a derivation
on $UL$, making $U(L,\del)$ into a \dga.
\subsection{Divided powers algebras}

Divided powers algebras arise here as the duals of universal
enveloping algebras.

\begin{defn}\label{def:dp-alg}
A \emph{divided powers algebra}, or $\Gamma$-algebra,
is a commutative graded algebra $A$, satisfying either
$A = A^{\geq 0}$ or $A = A^{\leq 0}$, equipped with set maps
\(
  \gamma^{k}:A^{2n} \rightarrow A^{2nk}
\)
for $k \geq 0$ and $n \neq 0$
satisfying the following list of conditions.
\begin{enumerate}
\item \label{dp1} $\gamma^{0}(a) = 1; \gamma^{1}(a)=a$
  for $a \in A$;
\item \label{dp2} 
  $\gamma^{k}(a+b) = \sum_{j=0}^{k}\gamma^{j}(a)\gamma^{k-j}(b)$
  for $a, b \in A^{2n}$;
\item \label{dp3} $\gamma^{j}(a)\gamma^{k}(a) =
  \displaystyle{\binom{j+k}{j}}\gamma^{j+k}(a)$ for $a \in A^{2n}$;
\item \label{dp4} $\gamma^{j}(\gamma^{k}(a)) =
  \displaystyle{\frac{(jk)!}{j!^{k}k!}}\gamma^{j+k}(a)$
  for $a \in A^{2n}$;
\item \label{dp5} $\gamma^{k}(ab) = \left\{
  \begin{array}{ll}
    a^{k}\gamma^{k}(b) & \mbox{if $|a|$ and $|b|$ even, $|b| \neq 0$,} \\
    0 & \mbox{if $|a|$ and $|b|$ odd.}
  \end{array}\right.$
\end{enumerate}
\end{defn}

A $\Gamma$-\emph{morphism} is an algebra morphism which respects the
divided powers operations.
A $\Gamma$-\emph{derivation} on a $\Gamma$-algebra $A$ is a derivation 
$\theta$ on $A$ satisfying 
$\theta(\gamma^{k}(a)) = \theta(a)\gamma^{k-1}(a)$
for $a \in A^{2n}$, $k \geq 1$.
A differential graded $\Gamma$-algebra, or \dgg,
is a pair $(A,\del)$, where $A$ is a $\Gamma$-algebra,
and $\del$ is a $\Gamma$-derivation of degree $-1$ satisfying
$\del\del = 0$.

Let $V$ be a free graded $R$-module.
Let $\Gamma^{k}(V)$ be the graded submodule of $T_{C}^{k}V$ 
of elements fixed by the action of the symmetric group $S_{k}$.
Then $\Gamma(V)=\dirsum_{k}\Gamma^{k}(V)$ is a Hopf subalgebra
of $T_{C}(V)$, called the \emph{free} $\Gamma$-\emph{algebra} on $V$.
Divided powers are defined on $\Gamma(V)$
by
\begin{enumerate}
  \item $\gamma^{0}(v) = 1$, $\gamma^{1}(v)=v$ for $v \in V$,
  \item $\gamma^{k}(v) = \underbrace{[v|\cdots|v]}_{k \mbox{\ times}}$
    for $v \in V^{2n}$
\end{enumerate}
and then extending via conditions (\ref{dp4}) and (\ref{dp5}) of
Definition~\ref{def:dp-alg}.
If $f:V \rightarrow A$ is any linear
map of degree zero from $V$ into a $\Gamma$-algebra $A$,
then $f$ extends to a unique $\Gamma$-morphism
$\bar{f}:\Gamma(V) \rightarrow A$.
If $V$ is $R$-free on a countable, well-ordered basis 
$\{ v_{i} \}$, then $\Gamma(V)$ is $R$-free, with basis
consisting of elements 
\(
  \gamma^{k_{1}}(v_{1})\cdots\gamma^{k_{s}}(v_{s})
\)
where $k_{j} \geq 0$ and $k_{j} = 0$ or $1$ if $|v_{j}|$ is odd.

If $V \tensor W \stackrel{\langle,\rangle}{\rightarrow} R$ is
a pairing, then there is an induced pairing
\begin{equation}\label{tensorpairing}
  TV \tensor T_{C}W \rightarrow R
\end{equation}
given by $\langle T^{j}V,T_{C}^{k}W \rangle = 0$ if $j \neq k$,
and
\begin{equation}\label{dual-rule}
  \langle v_{1} \tensor \cdots \tensor v_{k},
    [w_{1}|\cdots|w_{k}] \rangle
  = \pm \langle v_{1},w_{1} \rangle \cdots 
      \langle v_{k},w_{k} \rangle
\end{equation}
where $\pm$ is the sign of the permutation
\[
  v_{1},\ldots,v_{k},w_{1},\ldots,w_{k}
  \mapsto v_{1},w_{1},\ldots,v_{k},w_{k}.
\]
The pairing (\ref{tensorpairing}) in turn induces a pairing
\begin{equation}\label{lambda-gamma}
  \Lambda V \tensor \Gamma W \rightarrow R.
\end{equation}
Suppose that $V$ is $R$-free of finite type, 
$V = V_{<0}$ or $V = V_{>0}$, and $W = V^{\dual}$.
Then (\ref{tensorpairing})
and (\ref{lambda-gamma}) induce Hopf algebra isomorphisms
\(
  T_{C}(V^{\dual}) \iso (TV)^{\dual}
\)
and
\(
  \Gamma(V^{\dual}) \iso (\Lambda V)^{\dual}.
\)
\subsection{The Cartan--Chevalley--Eilenberg--Cartan complex}

Denote by $B(A)$ the bar construction on the augmented
\dga\ $(A,\del)$ (\cite{halperin:92}, Section 1);
recall that the underlying coalgebra of $B(A)$ is
$T_{C}(s\bar{A})$, where $\bar{A}$ is the augmentation ideal.
Let $(L,\del)$ be a \dgl.
Then $\Gamma(sL)\subset\Gamma(s\overline{UL})\subset B(UL)$
and $(\Gamma(sL),\del_{0}+\del_{1})$ is a sub-\dgc\ of
$B(UL)$, denoted by $C_{*}(L,\del)$, called the
\emph{chains on} $(L,\del)$.

The \emph{Cartan--Chevalley--Eilenberg--Cartan complex}
on $(L,\del)$ is the commutative cochain algebra
$C^{*}(L,\del) = (\Lambda V,d)$,
dual to $C_{*}(L,\del)$,
where $V=(sL)^{\dual}$, and the differential $d$ is the
sum of derivations $d_{0}$ and $d_{1}$.
The \emph{linear part} $d_{0}$ preserves word length and
is dual to $\del$ in that
\(
  \langle d_{0}v,sx \rangle = (-1)^{|v|}\langle v, s\del x \rangle
\)
for $v \in V$, $x \in L$.
The \emph{quadratic part} $d_{1}$ increases word length by one
and is dual to the Lie bracket in $L$:
\begin{equation}
  \langle d_{1}v,sx\cdot sy \rangle 
    = (-1)^{|sy|} \langle v, s[x,y] \rangle  \label{eq:quad}
\end{equation}
where the pairing is (\ref{lambda-gamma}) above
with $W = sL = V^{\dual}$.
We will usually refer to the Cartan--Chevalley--Eilenberg--Cartan
complex as the \emph{cochains on} $(L,\del)$.
\subsection{Bockstein spectral sequences}

Fix a prime $p$.
Let $C$ be a free chain complex over $\zlocp$.
Applying $C \tensor -$ to the short exact sequence of coefficient modules
\[
  0 \rightarrow \zlocp \stackrel{\times p}{\rightarrow}
    \zlocp \rightarrow \fp \rightarrow 0
\]
leads to a long exact sequence in homology which may be wrapped into the
exact couple
\[
  \begin{diagram}
  \node{H_{*}(C)} \arrow[2]{e} \node[2]{H_{*}(C)} \arrow{sw} \\
  \node[2]{H_{*}(C;\fp)}\arrow{nw}
  \end{diagram}
\]
from which we get the \emph{homology Bockstein
spectral sequence modulo} $p$ of $C$,
$(E^{r}(C),\beta^{r})$,
mod $p$ \bss\ for short~\cite{browder:61}.
If $C = C_{*}(X)$ is the normalized singular chain complex of a
space $X$, then we refer to the homology \bss\ mod $p$ of
$C_{*}(X)$ as the mod $p$ homology \bss\ of $X$, denoted 
$(E^{r}(X),\beta^{r})$.

There is the corresponding notion of 
\emph{cohomology Bockstein spectral sequence}
defined in the obvious manner, using the functor
$\Hom(C,-)$ rather than $C \tensor -$.

The mod $p$ \bss\ of $C$ measures $p$-torsion in $H_{*}(C)$:
if $x,y\in E^{r}$, $x \neq 0$, satisfy $\beta^{r}(y)=x$, then $x$ represents
a torsion element of order $p^{r}$ in $H_{*}(C)$.

\emph{Notation.}  
If $c\in C$ is such that $[\bar{c}]\in E^{1}$ lives until
the $E^{r}$ term then we will denote the corresponding element of 
$E^{r}$ by $[c]_{r}$.
\subsection{Acyclic closures and minimal models}

(Reference:~\cite{halperin:92},  Sections 2 and 7.)
Consider the graded algebra $\Lambda V \tensor \Gamma(sV)$
over $R$.
Extend the divided powers operations on $\Gamma(sV)$ to
$R \dirsum \Lambda V \tensor \Gamma^{+}(sV)$ via rule 5 of 
Definition~\ref{def:dp-alg}.

\begin{defn}(\cite{halperin:92}, Section 2)
An acyclic closure of the \dga\ $(\Lambda V,d)$ is a 
\dga\ of the form $C = (\Lambda V \tensor \Gamma(sV),D)$
in which $D$ is a $\Gamma$-derivation restricting to $d$ in
$\Lambda V$ and $H(C) = H^{0}(C) = R$.
\end{defn}
 
Let $(L,\del)$ be a connected \dgl\ over $R$ which is $R$-free of finite
type.
Then $C^{*}(L) = (\Lambda V,d)$ where $V = (sL)^{\dual}$.
Let $C$ be an acyclic closure for $C^{*}(L)$,
and set $(\Gamma(sV),\bar{D}) = R \tensor_{C^{*}(L)}C$.
By the work of Halperin in~\cite{halperin:92},
we identify $H(UL) = H([\Gamma(sV),\bar{D}]^{\dual})$
and $UL = (\Gamma(sV),\bar{D})^{\dual}$.

Let $R = \zlocp$ or $R = \fp$, and
consider a commutative algebra of the form $(\Lambda W,d)$
over $R$,
where $W = W^{\geq 2}$ is $R$-free and of finite type.
We may write the differential as a sum
$d = \sum_{j \geq 0}d_{j}$ where $d_{j}$ raises wordlength
by $j$.

\begin{defn}
If $R = \zlocp$, the \dga\ $(\Lambda W,d)$ above is 
$\zlocp$-\emph{minimal} if $d_{0}:W \rightarrow pW$.
If $R = \fp$, $(\Lambda W,d)$ is $\fp$-\emph{minimal}
if $d_{0} = 0$.
\end{defn}

Suppose $(A,\del)$ is a cochain algebra satisfying
$H^{0}(A) = R$, $H^{1}(A)=0$, $H^{2}(A)$ is $R$-free,
and $H^{*}(A)$ is of finite type.
Then by~\cite{halperin:92}, Theorem 7.1, there exists a quasi-isomorphism
$m:(\Lambda W,d)\quism(A,\del)$ from an $R$-minimal algebra.
This quasi-isomorphism is called a \emph{minimal model}.

Associated to an $\fp$-minimal model
$m:(\Lambda W,d)\quism(A,\del)$
is its \emph{homotopy Lie algebra}, $E$.
As a graded vector space, $E = (sW)^{\dual}$;
the bracket is defined by the relation
\[
  \langle w, s[x,y] \rangle 
    = (-1)^{|sy|} \langle d_{1}w, sx\cdot sy \rangle
\]
for $w \in W$, $x,y \in E$. 


\section{The image of $H(L)\rightarrow H(UL)$}\label{s:lie-morph}

Let $(L,\del)$ be a connected \dgl\ over $\fp$ of finite type.
By~\cite{halperin:92}, the choice
of minimal model $m:(\Lambda W,d)\quism C^{*}(L)$
determines an isomorphism of graded Hopf algebras,
$H(UL) \iso UE$, where $E$ is the 
homotopy Lie algebra of $m$.

\begin{prop}\label{prop:theta}
With the notation above, the image of 
$H(\iota):H(L)\rightarrow H(UL)$ lies in $E$.
\end{prop}

\begin{proof}
It suffices to prove that the following diagram commutes.
\begin{equation}
\begin{diagram}\label{d:theta}
\node{H(L,\del)}\arrow{e,t}{H(\iota_{L})}\arrow{s,l}{\theta}
  \node{H(UL)}\arrow{s,r}{\iso} \\
\node{E}\arrow{e,t}{\iota_{E}}\node{UE}
\end{diagram}
\end{equation}

Recall that $C^{*}(L,\del) = (\Lambda V,d)$,
where $V = (sL)^{\dual}$ and $d = d_{0} + d_{1}$.
Recall further that the minimality condition on $(\Lambda W,d)$
implies that the linear part of its differential vanishes.
The \emph{linear part} of $m$ is the linear map
$m_{0}:(W,0)\rightarrow (V,d_{0})$
defined by the condition $m - m_{0}:W \rightarrow \Lambda^{\geq 2}V$.
Recall that $E = (sW)^{\dual}$ and
$UE = \Gamma(sW)^{\dual}$ (\cite{halperin:92}, Theorem 6.2).

The model $m$ extends to a morphism of constructible acyclic closures
(\cite{halperin:92}, Section 2)
\(
  \hat{m}:(\Lambda W\tensor \Gamma(sW),D)
    \rightarrow (\Lambda V \tensor \Gamma(sV),D)
\)
by Proposition 2.7 of~\cite{halperin:92}.
Since $(\Lambda W,d)$ is $\fp$-minimal,
$d_{0}=0$.
By Corollary 2.6 of~\cite{halperin:92},
$d_{0}=0$ is equivalent to $\bar{D}=0$
in $(\Gamma(sW),\bar{D})$.
Apply $\fp\tensor_{m}\!-$ to $\hat{m}$
to get a $\Gamma$-morphism
$\bar{m}:(\Gamma(sW),0)\rightarrow(\Gamma(sV),\bar{D})$.

Let 
$\pi_{L}:(\Gamma(sV),\bar{D})\twoheadrightarrow s(V,d_{0})$
and
$\pi_{E}:(\Gamma(sW),0)\twoheadrightarrow s(W,0)$
be the projections.
The maps $\pi_{L}$ and $\pi_{E}$ fit into the diagram
\begin{equation}
\begin{diagram}\label{d:pi}
  \node{(\Gamma(sW),0)}\arrow{e,t}{\pi_{E}}
    \arrow{s,l}{\bar{m}}\node{s(W,0)}
    \arrow{s,r}{sm_{0}}      \\
  \node{(\Gamma(sV),\bar{D})}\arrow{e,t}{\pi_{L}}
    \node{s(V,d_{0})}
\end{diagram}
\end{equation}
For $w \in W$, Proposition 2.7 of~\cite{halperin:92}
states that
$\hat{m}(1 \tensor sw) - 1 \tensor sm_{0}w$
has total wordlength at least two.
It follows that $\bar{m}(sw) - sm_{0}w$ has
$\Gamma(sV)$-wordlength at least two, so
$\pi_{L}(\bar{m}(sw)) = sm_{0}w = sm_{0}(\pi_{E}(sw))$,
so Diagram (\ref{d:pi}) commutes.
Dualize and pass to homology to get (\ref{d:theta}).
\end{proof}


\section{Bockstein spectral sequence of a 
universal enveloping algebra}\label{s:bss}

In this section, we prove the main algebraic result,
Theorem~\ref{thm:main}.
Unless otherwise stated, our ground ring will be
$\zlocp$, the integers localized at $p$.

Let $(\Lambda W,d)$ be a minimal Sullivan algebra over $\zlocp$.
Let $C = (\Lambda W \tensor \Gamma(sW),D)$ be a constructible acyclic
closure for $(\Lambda W,d)$ (\cite{halperin:92}, Section 2).
Let $(\Gamma(sW),\bar{D})$ be the quotient
$\zlocp\tensor_{(\Lambda W,d)}C$.
$C \tensor \fp$ is a constructible acyclic closure for
$(\Lambda W,d)\tensor\fp$.
Since $(\Lambda W,d)$ is $\zlocp$-minimal, $p|d_{0}$,
so the linear part of the differential vanishes
in $(\Lambda W,d)\tensor \fp$.
It follows by Corollary 2.6 of~\cite{halperin:92} that the
differential in $(\Gamma(sW),\bar{D})\tensor\fp$ is null,
so that $p | \bar{D}$.
Set $E^{r} = E^{r}([\Gamma(sW),\bar{D}]^{\dual})$ and 
$E_{r} = E_{r}(\Gamma(sW),\bar{D})$.
Let $\rho:\Gamma(sW)\rightarrow \Gamma(sW)\tensor\fp = E_{1}$
be the reduction homomorphism.

\begin{prop}\label{prop:main}
With the hypotheses and notation above, for $r\geq 1$, 
the following statements hold.

\begin{enumerate}
\item $E_{r}$ is isomorphic to a free divided powers
algebra, 
\item there is a $\Gamma$-morphism
$g_{r}:E_{r}\rightarrow E_{1}$ such that if $g(z)=\rho(a)$
for some $z \in E_{r}$, $a \in \Gamma(sW)$, then $z = [a]_{r}$.
\item there is a graded Lie algebra $L^{r}$ such that 
$(E^{r},\beta^{r}) = U(L^{r},\beta^{r})$.
\end{enumerate}
\end{prop}

\begin{lem}\label{lem:key}
Let $(UL,\del)$ be a \dgh\ over $\fp$ of finite type,
so
$(UL)^{\dual} = \Gamma V$ as an algebra.
If $\del^{\dual}$
is a $\Gamma$-derivation, then $\del(L)\subset L$.
\end{lem}

\begin{proof}[Proof of Lemma~\ref{lem:key}]
It suffices to prove the dual statement, namely that
$\del^{\dual}:\Gamma V \rightarrow \Gamma V$ factors over
the surjection $\pi:\Gamma V \rightarrow V$ to induce a
differential in $V$.
But $\ker(\pi)$ consists of products along with elements
of the form $\gamma^{p^{k}}(v)$ for $v \in V$, $k \geq 1$.
Since $\del^{\dual}$ is a $\Gamma$-derivation,
$\del^{\dual}(\gamma^{p^{k}}(v)) = \del^{\dual}(v)\gamma^{p^{k}-1}(v)$
is a product.
It follows that $\del^{\dual}(\ker(\pi)) \subset \ker(\pi)$,
completing the proof.
\end{proof}

\begin{proof}[Proof of Proposition~\ref{prop:main}]
We proceed by induction.
For $r=1$, let $W_{1} = W \tensor \fp$.
Since $p|\bar{D}$, $E_{1} = \Gamma(sW_{1})$,
establishing the first statement.
For the second statement we may take $g_{1}$ to be the identity map.
Let $L^{1}$ be the homotopy Lie algebra of 
$(\Lambda W,d)\tensor\fp$.
Apply Theorems 6.2 and 6.3 of~\cite{halperin:92} to
the minimal algebra $(\Lambda W,d)\tensor\fp$
to get a graded Hopf algebra isomorphism
$E^{1} = (\Gamma(sW_{1}))^{\dual} \iso UL^{1}$.

Since $p|\bar{D}$ in $(\Gamma(sW),\bar{D})$,
$\beta_{1} = \bar{D}/p$ (reduced modulo $p$).
Thus because $\bar{D}$ is a $\Gamma$-derivation,
so is $\beta_{1}$.
Since $\Gamma(sW_{1})$ is the $\Gamma$-algebra dual to
$UL^{1}$ it follows by Lemma~\ref{lem:key} 
that $\beta^{1}:L^{1}\rightarrow L^{1}$
and so $E^{1} = U(L^{1},\beta^{1})$.

Now suppose the three statements are established for $r-1$.
Let $C(r-1)=(\Lambda W_{r-1}\tensor\Gamma(sW_{r-1}),D)$ 
be a constructible acyclic closure for 
$C^{*}(L^{r-1},\beta^{r-1})=(\Lambda W_{r-1},d)$.
By Lemma 5.4 of~\cite{halperin:92}, there is a chain isomorphism
\(
  \gamma_{r-1}:(E^{r-1},\beta^{r-1})= U(L^{r-1},\beta^{r-1})
    \stackrel{\iso}{\rightarrow} (\Gamma(sW_{r-1}),\bar{D})^{\dual}
\).
Fix a well-ordered basis $\{x_{j}\}$ of $L^{r-1}$;
this determines a dual basis $\{ sw_{j} \}$ of $sW_{r-1}$.
The isomorphism $\gamma_{r-1}$ identifies the 
Poincar\'e--Birkhoff--Witt basis element 
$x_{1}^{k_{1}}\cdots x_{j}^{k_{j}}$
of $UL^{r-1}$ as a
dual basis to the basis element
$\gamma^{k_{1}}(sw_{1})\cdots\gamma^{k_{j}}(sw_{j})$
of $\Gamma(sW_{r-1})$.
It follows that $\gamma_{r-1}$ factors as the composition
of \dgc\ isomorphisms 
$UL^{r}\stackrel{\iso}{\rightarrow}\Lambda(L^{r})^{\flat}
\stackrel{\iso}{\rightarrow}(\Gamma(sW_{r-1}))^{\dual}$.
Since $L^{r-1}$ and $W_{r-1}$ are finite type, we dualize
to obtain the \dga\ isomorphism
\(
  \alpha_{r-1}:(\Gamma(sW_{r-1}),\bar{D})
    \stackrel{\iso}{\rightarrow}(E_{r-1},\beta_{r-1})
\).

Let $m_{r}:(\Lambda W_{r},d)\quism C^{*}(L^{r-1},\beta^{r-1})$
be a minimal model.
Let $C'(r)$ be a constructible acyclic closure
of $(\Lambda W_{r},d)$ 
(\cite{halperin:92}, Section 2).
Since $(\Lambda W_{r},d)$ is $\fp$-minimal, $d_{0}=0$,
so by Corollary 2.6 of~\cite{halperin:92}, the differential
in $(\Gamma(sW_{r}),\bar{D}_{(r)})=\fp \tensor_{\Lambda W}C'(r)$ is zero.
By~\cite{halperin:92}, Proposition 2.7, $m_{r}$ induces a $\Gamma$-morphism
\(
  \bar{m}_{r}:(\Gamma(sW_{r}),0)\quism(\Gamma(sW_{r-1}),\bar{D}).
\)
Since $\fp$ is a field, by Lemma 3.3 of~\cite{halperin:92},
we may identify $H(\bar{m}_{r})$ with
$\tor^{m_{r}}(\fp,\fp)$,
where $\tor$ is the differential torsion functor~\cite{f-h-t:95}.
Therefore
since $m_{r}$ is a quasi-isomorphism,
$H(\bar{m}_{r})$ is an isomorphism. 
Let $\alpha_{r}:\Gamma(sW_{r})\stackrel{\iso}{\rightarrow}E_{r}$
be the composition of algebra isomorphisms
\[
  \Gamma(sW_{r})\xrightarrow{H(\bar{m}_{r})}
  H(\Gamma(sW_{r-1}),\beta_{r-1})
  \xrightarrow{H(\alpha_{r-1})}E_{r}
\]
to establish the first statement.
Note that $\alpha_{r}$ will be the isomorphism dual to $\gamma_{r}$
in the next stage of the induction.

Setting $f = \alpha_{r-1}\bar{m}_{r}\alpha_{r}^{-1}$,
we get the commutative diagram
\[
  \begin{diagram}
  \node{(\Gamma(sW_{r}),0)}\arrow{e,t}{\bar{m}_{r}}
    \arrow{s,lr}{\alpha_{r}}{\iso}
    \node{(\Gamma(sW_{r-1}),\bar{D})}
    \arrow{s,lr}{\alpha_{r-1}}{\iso}   \\
  \node{(E_{r},0)}\arrow{e,t}{f}\node{(E_{r-1},\beta_{r-1})}
  \end{diagram}
\]
and it follows from the definitions that $H(f)$ is the identity
on $E_{r}$.

By the inductive hypothesis, there exists a $\Gamma$-morphism
$g_{r-1}:E_{r-1} \rightarrow E_{1}$ such that $z = [a]_{r-1}$
whenever $z \in E_{r-1}$, $a \in \Gamma(sW)$ satisfy $g(z)=\rho(a)$.
We now show that $g_{r}:=g_{r-1}f$ satisfies statement 2.
For $u \in E^{r}$ choose $a \in \Gamma(sW)$ so that
$g_{r-1}(f(u))=\rho(a)$.
Then $f(u)=[a]_{r-1}$, hence $\beta_{r-1}[a]_{r-1}=0$ and
$[a]_{r}\in E^{r}$ is defined.
Since $f$ induces the identity in homology,
$f([a]_{r}) = [a]_{r-1} + \beta_{r-1}(v)$
for some $v \in E^{r-1}$.
Thus $f(u - [a]_{r})=\beta_{r-1}(v)$,
so $u - [a]_{r}$ is a boundary in $(E_{r},0)$, whence $u=[a]_{r}$.
This establishes the second statement.

The model $m_{r}$ determines an isomorphism
$E^{r} = H(U(L^{r-1},\beta^{r-1}))
\iso UL^{r}$, where $L^{r}$ is the homotopy Lie algebra for the
model $m_{r}$.
As a graded vector space, $L^{r} = (sW_{r})^{\dual}$.

Let $u \in E_{r}$, and suppose for some 
$a \in \Gamma(sW)$ that $\rho(a) = g_{r}(u)$.
Then $u = [a]_{r}$,
so $\bar{D}a=p^{r}b$ for some $b \in \Gamma(sW)$.
Thus
$\beta_{r}(u) = [b]_{r}$.
Since $g_{r}$ and $\rho$ are $\Gamma$-morphisms,
$\rho(\gamma^{j}(a))=g_{r}(\gamma^{j}(u))$ so
$\gamma^{j}(u) = [\gamma^{j}(a)]_{r}$.
Furthermore, $\bar{D}(\gamma^{k}(a)) = p^{r}b\cdot\gamma^{k-1}(a)$
so
\[
  \beta_{r}\gamma^{k}(u)
  = \beta_{r}[\gamma^{k}(a)]_{r}
  = [ b \cdot \gamma^{k-1}(a) ]_{r}
  = [b]_{r}[\gamma^{k-1}(a)]_{r}
  = \beta_{r}(u)\cdot\gamma^{k-1}(u).
\]
By Lemma~\ref{lem:key}, this establishes the third statement, 
completing the inductive
step and the proof.
\end{proof}

\begin{proof}[Proof of Theorem~\ref{thm:main}]
Let $m:(\Lambda W,d) \quism C^{*}(L,\del)$ be a minimal model.
Recall that the underlying algebra of $C^{*}(L,\del)$ is $\Lambda V$,
where $V = (sL)^{\dual}$.
Let $(\Lambda W \tensor \Gamma(sW),D)$ and 
$(\Lambda V \tensor \Gamma(sV),D)$  be constructible acyclic
closures for $(\Lambda W,d)$ and $C^{*}(L,\del)$, respectively.
The model $m$ determines a $\Gamma$-morphism
$\bar{m}:(\Gamma(sW),\bar{D})\rightarrow(\Gamma(sV),\bar{D})$
where $H(\bar{m}^{\dual})$ is an isomorphism.
The composition
\[
  U(L,\del) \stackrel{\iso}{\rightarrow} 
    (\Gamma(sV),\bar{D})^{\dual} \quism
    (\Gamma(sW),\bar{D})^{\dual}
\]
induces an isomorphism of Bockstein spectral sequences,
establishing the first statement.

The reduced minimal model
$m \tensor \fp:(\Lambda W,d)\tensor\fp\quism C^{*}(L,\del)\tensor\fp$
has homotopy Lie algebra $L^{1}$, so by Proposition~\ref{prop:theta}, 
$\image{E^{1}(\iota)}\subset L^{1}$.
Suppose that $\image{E^{r-1}(\iota)}\subset L^{r-1}$.
Let $\iota^{(r-1)}:L^{r-1}\hookrightarrow UL^{r-1}$ be the
inclusion.
Then $\image{E^{r}(\iota)}\subset\image{H(\iota^{(r-1)})}$.
The homotopy Lie algebra of the minimal model
$m_{r}:(\Lambda W_{r},d)\quism C^{*}(L^{r-1},\beta^{r-1})$
is $L^{r}$, so Proposition~\ref{prop:theta} states that
$\image{H(\iota^{(r-1)})}\subset L^{r}$, completing the induction
and the proof.
\end{proof}

\begin{proof}[Proof of Theorems~\ref{thm:loopbss} 
  and~\ref{thm:transf}]
Anick in~\cite{anick:89} proves that there is a \dgl\ $L_{X}$
and a \dgh\ quasi-isomorphism 
$UL_{X}\rightarrow C_{*}(\Omega X;\zlocp)$.
Thus as Hopf algebras, for $r \geq 1$,
$E^{r}(UL_{X}) = E^{r}(\Omega X)$ and 
$E_{r}(UL_{X}) = E_{r}(\Omega X)$.
The result follows by applying Theorem~\ref{thm:main} to
the \dgl\ $L_{X}$.
\end{proof}

%
%
\section{Morphisms of universal enveloping algebras}\label{s:morphisms}

Let $R$ be a commutative ring, $L_{1}$ and $L_{2}$ connected,
graded Lie algebras over $R$ which are $R$-free of finite type.
Let $\varphi:UL_{1} \rightarrow UL_{2}$ be a morphism of
Hopf algebras.
The purpose of this section is to prove
\begin{prop}
$\varphi(L_{1})\subset L_{2}$ if and only if
$\varphi^{\dual}:(UL_{2})^{\dual}\rightarrow(UL_{1})^{\dual}$
is a $\Gamma$-morphism.
\end{prop}

\begin{proof}
Observe that, for $j = 1,2$, the sequence of functors
\[
  L_{j} \rightsquigarrow C^{*}(L_{j})
    \rightsquigarrow \Gamma V_{j} 
      = R \tensor_{C^{*}(L_{j})} (C^{*}(L_{j}) \tensor \Gamma V_{j})
    \rightsquigarrow (\Gamma V_{j})^{\dual} = UL_{j}
\]
identifies $UL_{j}$ as the natural dual of the free $\Gamma$-algebra
$\Gamma V_{j}$, with $L_{j}$ naturally dual to $V_{j}$, and
the inclusion $L_{j} \hookrightarrow UL_{j}$
naturally dual to the projection 
$\Gamma V_{j} \stackrel{\pi_{j}}{\twoheadrightarrow}V_{j}$.
Therefore if $\varphi|_{L_{1}}:L_{1}\rightarrow L_{2}$,
then $\varphi^{\dual}$ is a $\Gamma$-morphism.
Conversely, if $\varphi^{\dual}$ is a $\Gamma$-morphism,
then $\varphi^{\dual}(\ker \pi_{2}) \subset \ker\pi_{1}$,
so $\varphi(L_{1})\subset L_{2}$.
\end{proof}


\section{Examples}\label{s:examples}

We begin with a proposition to be used in both examples.

\begin{prop}\label{prop:example}
Define a \dgl\ over $\fp$ by
$(L, \del) = (\ab(e,f), \del f = e)$,
where $|f| = 2n$.
Then
$C^{*}(L,\del)=(\Lambda(x,y),d)$ with $dx=y$
and $|x| = 2n$.
A minimal model
$m:(\Lambda(x_{1},y_{1}),0)\quism C^{*}(L,\del)$,
given by
$x_{1} \mapsto x^{p}$
and
$y_{1} \mapsto x^{p-1}y$,
induces isomorphisms
\(
  \Gamma(sx_{1},sy_{1})\stackrel{\iso}{\rightarrow}
    H([UL]^{\dual})
\)
and
\(
   H(UL)\stackrel{\iso}{\rightarrow}
   U \ab(e_{1},f_{1})
\)
with $|e_{1}|=|sx_{1}|= 2np - 1$, $|f_{1}|=|sy_{1}|=2np$.
\end{prop}

\begin{proof}
Straightforward. 
\end{proof}

\begin{ex}
Let $L = \ab(e,f)$
over $\zlocp$ on
generators $e$ and $f$ of degrees $2n-1$ and $2n$, respectively.
Set $\del f = pe$.
Applying Proposition~\ref{prop:example} recursively, we have
$E_{r}(UL) = \Gamma(sx_{r},sy_{r})$ and
$E^{r}(UL) = U \ab(e_{r},f_{r})$,
with
$|e_{r}| = |sx_{r}| = 2np^{r}-1$,
$|f_{r}| = |sy_{r}| = 2np^{r}$,
$\beta_{r}(sx_{r}) = sy_{r}$,
and
$\beta^{r}(f_{r}) = e_{r}$,
while the sequence $E^{r}(L)$ collapses after the first term.
\end{ex}

\begin{ex}
Define a \dgl\ $(L,\del)$ over $\zlocp$ by
$L = \ab(e,f,g)$, where $|e| = 2n-1$,
$|f| = |g|= 2$, and $\del(f) = pe$.
Then $L^{1} = \ab(e,f,g)$ (over $\fp$),
with $\beta^{1}(f)=e$,
and $C^{*}(L^{1},\beta^{1})=(\Lambda(x,y),dx=y)\tensor(\Lambda(z),0)$.
Recall the model $m$ from Proposition~\ref{prop:example}.
Define \dga\ morphisms 
$i,j:(\Lambda(z),0)\rightarrow C^{*}(L^{1},\beta^{1})$
by $i(z)=z$, $j(z)=z+y$.
Then $\varphi = m \tensor i$ and
$\psi = m \tensor j$ are minimal models,
both with homotopy Lie algebra
$L^{2} = \ab(a,b,c)$,
$|a|=2np-1$,$|b|=2np$, and $|c| = 2n$.
The two models determine Hopf algebra isomorphisms
$\varphi^{*},\psi^{*}:H(UL^{1})\rightarrow UL^{2}$,
given by
$\varphi^{*}[ef^{p-1}] = \psi^{*}[ef^{p-1}] = a$,
$\varphi^{*}[g] = \psi^{*}[g] = c$,   
$\varphi^{*}[f^{p}]=b$, and
$\psi^{*}[f^{p}]=b+c^{p}$.
The algebra isomorphism 
$\psi^{*}(\varphi^{*})^{-1}:U\ab(a,b,c) \rightarrow U\ab(a,b,c)$ 
is not of the
form $U\theta$ for any Lie algebra morphism
$\theta: \ab(a,b,c) \rightarrow \ab(a,b,c)$.
Therefore the construction involved in
Theorem~\ref{thm:main} is not natural.
\end{ex}

\bibliographystyle{amsplain}
\bibliography{loopy}

\end{document}